\documentclass[12pt]{article}
\usepackage{amsmath,amssymb,amsthm, comment}

\newtheorem{theorem}{Theorem}

\newtheorem{corollary}[theorem]{Corollary}

\def\barr{\begin{array}}
\def\earr{\end{array}}

\title{A result on the number of cyclic subgroups of a finite group}
\author{Marius T\u arn\u auceanu}
\date{March 13, 2020}

\begin{document}

\maketitle

\begin{abstract}
Let $G$ be a finite group and $\alpha(G)=\frac{|C(G)|}{|G|}$\,, where $C(G)$ denotes the set of cyclic subgroups of $G$.
In this short note, we prove that $\alpha(G)\leq\alpha(Z(G))$ and we describe the groups $G$ for which the equality occurs.
This gives some sufficient conditions for a finite group to be $4$-abelian or abelian.
\end{abstract}
\vspace{1mm}
{\small
\noindent
{\bf MSC2000\,:} Primary 20D60; Secondary 20D15, 20F18.

\noindent
{\bf Key words\,:} finite groups, $p$-groups, number of cyclic subgroups.}

\section{Introduction}
Let $G$ be a finite group, $C(G)$ be the set of cyclic subgroups of $G$ and $Z(G)$ be the center of $G$.
In 2018, M. Garonzi and I. Lima introduced in their paper \cite{1} the function
\begin{equation}
\alpha(G)=\frac{|C(G)|}{|G|}\,.\nonumber
\end{equation}Since then many authors have studied the properties of this function and its relations with the structure of $G$.
Note that we have
\begin{equation}
\alpha(G)=\frac{1}{|G|}\sum\limits_{x\in G}\frac{1}{\varphi(o(x))}\,,\nonumber
\end{equation}where $\varphi$ is Euler's totient function and $o(x)$ is the order of $x\in G$, showing that in fact
$\alpha(G)$ depends only on the element orders of $G$. Denote by $o(G)$ the average order of $G$, that is\newpage
\begin{equation}
o(G)=\frac{1}{|G|}\sum\limits_{x\in G}o(x).\nonumber
\end{equation}By Lemma 2.7 of \cite{4} (see also Corollary 2.6 of \cite{6}), this satisfies the following beautiful inequality
\begin{equation}
o(G)\geq o(Z(G)).\nonumber
\end{equation}

Our main result shows that a reversed inequality holds for the function $\alpha$ and gives a description of finite groups $G$ for which the equality occurs.

\begin{theorem}
Let $G$ be a finite group. Then
\begin{equation}
\alpha(G)\leq\alpha(Z(G)),\nonumber
\end{equation}and we have equality if and only if $G\cong G_1\times G_2$, where $G_1$ is a $2$-group with $G_1=\Omega_{\lbrace 1\rbrace}(G_1)Z(G_1)$ and $G_2$ is an abelian group of odd order.
\end{theorem}

The above condition $G_1=\Omega_{\lbrace 1\rbrace}(G_1)Z(G_1)$ means that $G_1$ has a transversal modulo $Z(G_1)$ consisting only of elements of order $1$ or $2$. Important examples of such $2$-groups are abelian $2$-groups and almost extraspecial $2$-groups (note that if $G_1$ is an almost extraspecial $2$-group, then $\alpha(G_1)=\alpha(Z(G_1)))=\frac{3}{4}$ by Theorem 4 of \cite{5}).

Since a finite group $G$ with $\alpha(G)=\alpha(Z(G))$ satisfies $\exp(G/Z(G))=2$, we infer the following corollary:

\begin{corollary}
Every finite group $G$ with $\alpha(G)=\alpha(Z(G))$ is $2$-central.
\end{corollary}

Also, Theorem 1 of \cite{2} implies that: 

\begin{corollary}
Every finite group $G$ with $\alpha(G)=\alpha(Z(G))$ is $4$-abelian. Moreover, if $|G|$ is odd, then $G$ is abelian.
\end{corollary}

Most of our notation is standard and will usually not be repeated here. Elementary notions and results on groups can be found in \cite{3}.

\section{Proof of Theorem 1}

Let $G/Z(G)=\{a_1Z(G)=Z(G),...,a_mZ(G)\}$, where $m=[G:Z(G)]$. Then
\begin{equation}
|C(G)|=\sum\limits_{x\in G}\frac{1}{\varphi(o(x))}=\sum_{i=1}^m\sum\limits_{x\in Z(G)}\frac{1}{\varphi(o(a_ix))}\,.
\end{equation}

We will prove that for every $i=2,...,m$ we have
\begin{equation}
\sum\limits_{x\in Z(G)}\frac{1}{\varphi(o(a_ix))}\leq\sum\limits_{x\in Z(G)}\frac{1}{\varphi(o(x))}\,.
\end{equation}Let $k_i=\min\{o(y)\mid y\in a_iZ(G)\}$ and $y_i\in a_iZ(G)$ such that $o(y_i)=k_i$. By the proof of Theorem A, we infer that
\begin{equation}
o(y_ix)=\frac{k_i}{(k_i,o(x))}\,o(x)\,\,\vdots\,\,o(x),\,\forall\, x\in Z(G).\nonumber
\end{equation}This leads to $\varphi(o(x))\mid\varphi(o(y_ix))$ and so $\varphi(o(x))\leq\varphi(o(y_ix))$, $\forall\, x\in Z(G)$. Since $a_iZ(G)=y_iZ(G)$, we obtain
\begin{equation}
\sum\limits_{x\in Z(G)}\frac{1}{\varphi(o(a_ix))}=\sum\limits_{x\in Z(G)}\frac{1}{\varphi(o(y_ix))}\leq\sum\limits_{x\in Z(G)}\frac{1}{\varphi(o(x))}\,,\nonumber
\end{equation}as desired.

Clearly, (1) and (2) imply
\begin{equation}
|C(G)|\leq m\!\!\!\sum\limits_{x\in Z(G)}\frac{1}{\varphi(o(x))}=m\,|C(Z(G))|,\nonumber
\end{equation}that is
\begin{equation}
\alpha(G)\leq\alpha(Z(G)).\nonumber
\end{equation}

Next, we remark that the equality $\alpha(G)=\alpha(Z(G))$ holds if and only if $$\varphi(o(y_ix))=\varphi(o(x)),\,\forall\,i=2,...,m,\,\forall\, x\in Z(G).$$By taking $x=1$, we get $\varphi(k_i)=1$, i.e. $k_i=2$. This shows that $G/Z(G)$ is an elementary abelian $2$-group. It follows that all Sylow $p$-subgroups of $G$ for $p$ odd are central and consequently $$G\cong G_1\times G_2,$$where $G_1$ is a $2$-group and $G_2$ is an abelian group of odd order. Obviously, we have $$Z(G)\cong Z(G_1)\times G_2 \mbox{ and } y_i\in\Omega_{\lbrace 1\rbrace}(G_1), \forall\, i=1,...,m,$$implying that $G_1=\Omega_{\lbrace 1\rbrace}(G_1)Z(G_1)$.

This completes the proof.$\qed$

\vspace*{3ex}\small

\hfill
\begin{minipage}[t]{5cm}
Marius T\u arn\u auceanu \\
Faculty of  Mathematics \\
``Al.I. Cuza'' University \\
Ia\c si, Romania \\
e-mail: {\tt tarnauc@uaic.ro}
\end{minipage}

\end{document}